\pdfoutput=1

\documentclass[12pt]{amsart}

\usepackage[indentafter]{titlesec}
\titleformat{name=\section}{}{}{0pt}{\scshape}

\usepackage[bookmarks=true,%
    colorlinks=true,%
    linkcolor=blue,%
    citecolor=blue,%
    filecolor=blue,%
    menucolor=blue,%
    urlcolor=blue,%
    breaklinks=true]{hyperref}
    
\usepackage{graphicx}

\usepackage[margin=2.5cm,letterpaper]{geometry}

\newcommand\myLangle{\langle~}
\newcommand\myRangle{~\rangle}

\renewcommand{\H}{\mathbb{H}}
\newcommand{\Q}{\mathbb{Q}}

\newcommand{\PSL}{\mathrm{PSL}}
\newcommand{\SL}{\mathrm{SL}}

\newcommand\mfdMarker[1]{\mathrm{\mathcal{#1}}}
\newcommand\MunivPrinCong[2]{\mfdMarker{U}^{#1}_{#2}}

\title{All Known Principal Congruence Links}
\author{M. D. Baker}
\author{M. Goerner}
\author{A. W. Reid}

\address{\newline
IRMAR,\newline
Universit\'e de Rennes 1,\newline
35042 Rennes Cedex,\newline
France.}
\email{mark.baker@univ-rennes1.fr }
\address{\newline 
Pixar Animation Studios, \newline
1200 Park Avenue,\newline
Emeryville, CA 94608, USA.}
\email{enischte@gmail.com}
\address{\newline
Department of Mathematics,\newline
Rice University,\newline
Houston, TX 77005, USA}
\email{alan.reid@rice.edu}

\begin{document}

\begin{abstract}
This report lists the link diagrams in $S^3$ for all principal congruence link complements for which such a link diagram is known. Several unpublished link diagrams are included. Related to this, we also include one link diagram for an arithmetic regular tessellation link complement.
\end{abstract}

\maketitle

\section{Introduction}

Given a square-free positive integer $d$ and an ideal $I$ in the ring of integers $O_d$ of the imaginary quadratic number field $\Q(\sqrt{-d})$, the principal congruence manifold (or orbifold) associated with $(d, I)$ is the quotient $\H^3/\Gamma(I)$ where $$\Gamma(I)=\ker\left(\PSL(2,O_d)\to \SL(2,O_d/I) / \{\pm 1\}\right).$$

There are only finitely principal congruence manifolds that are link complements in $S^3$ and the complete list was given in  \cite{bgr18:AllPrinCong}. We refer the reader to \cite{bakerReid:schwermer} and \cite{bakerReid:prinCong} for an introduction to and further background on congruence manifolds.

This report serves as a repository of link diagrams for the cases $(d,I)$ for which such a diagram is known. Since there are, in general, many links having the same complement, we show only one link diagram for a case. We plan to update this document as more principal congruence links become known. Please contact us if you have constructed a new link or know about a case not listed here. Furthermore, if you cite this article and refer to a particular link, do so by using the pair $(d,I)$ rather than the figure number as the figure numbers might potentially change as more principal congruence links become known.

All the link diagrams are also available files \cite[\href{http://unhyperbolic.org/prinCong/prinCong/Links/}{prinCong/Links/}]{goerner:data} which can be opened and viewed with SnapPy \cite{SnapPy}, e.g.,
\begin{center}
\begin{minipage}{12cm}
\begin{verbatim}
cd DIRECTORY_WITH_LINK_FILES
M = Manifold("pSL_0.5_plus_0.5_sqrt_minus_11.lnk")
M.plink()
\end{verbatim}
\end{minipage}
\end{center}

Thematically related are the regular tessellation link complements defined in \cite{goerner:regTessLinkComps} and we also include one link (Figure~\ref{fig:cubical}) even though it is not principal congruence.

A dot in a link diagram denotes a link component that intersects the paper plane orthogonally.

\clearpage

\section{$d=1$}

Missing: $(1, \myLangle 4 + \sqrt{-1}\myRangle)$.

\begin{figure}[h]
\centering
\begin{minipage}[t]{0.5\textwidth}
\begin{center}
\includegraphics[scale=0.5]{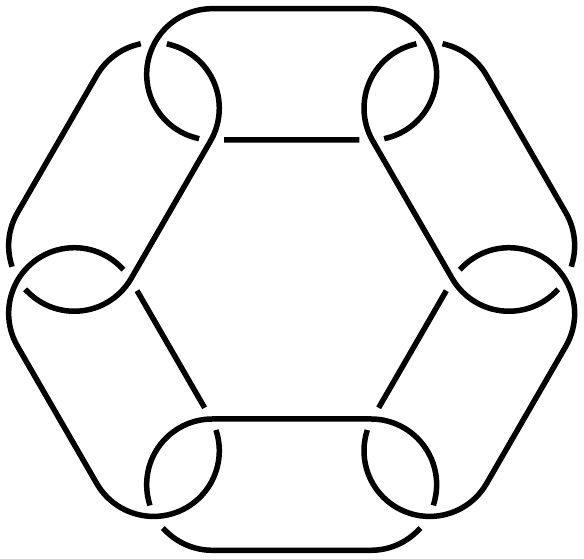}
\caption[caption]{$(1,\myLangle 2\myRangle)$.\\ \cite[Fig.~I-17]{baker:thesis}\\ \cite[Fig.~1a]{bakerReid:prinCong}.}
\end{center}
\end{minipage}%
\begin{minipage}[t]{0.5\textwidth}
\begin{center}
\includegraphics[scale=1.0]{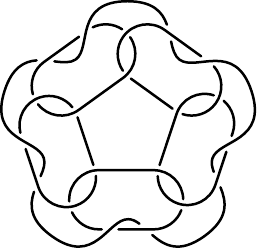}
\caption[caption]{$(1,\myLangle 2+\sqrt{-1}\myRangle)$.\\By second author.}
\end{center}
\end{minipage}
\end{figure}

\begin{figure}[h]
\begin{center}
\includegraphics[scale=1.0]{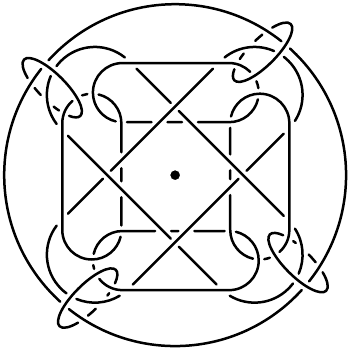}
\caption[caption]{$(1,\myLangle 2+2\sqrt{-1}\myRangle )$. By second author.}
\end{center}
\end{figure}

\begin{figure}[h]
\begin{center}
\includegraphics[scale=0.8]{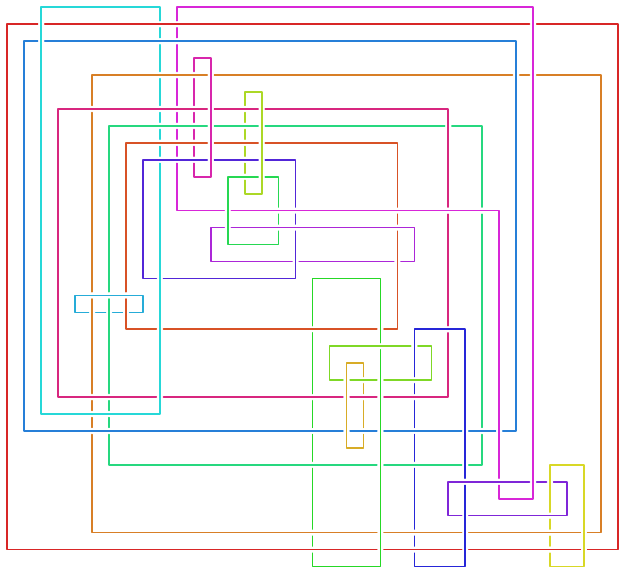}
\caption[caption]{$(1,\myLangle 3\myRangle )$. \cite{DOR:linkDiag}.}
\end{center}
\end{figure}

\begin{figure}[h]
\begin{center}
\includegraphics[scale=1.0]{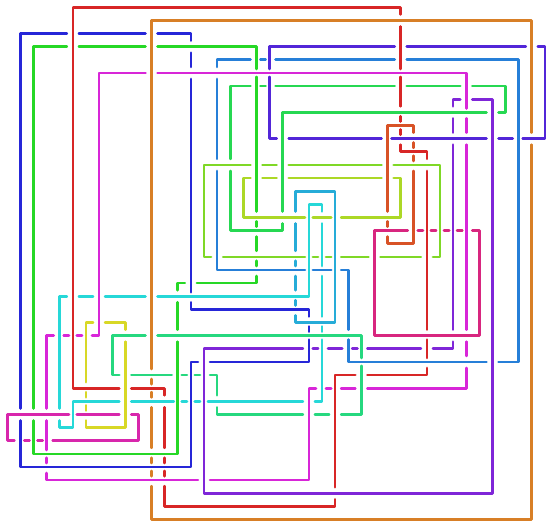}
\caption[caption]{$(1,\myLangle 3 + \sqrt{-1}\myRangle )$. \cite{DOR:linkDiag}.}
\end{center}
\end{figure}

\begin{figure}[h]
\begin{center}
\includegraphics[scale=1.5]{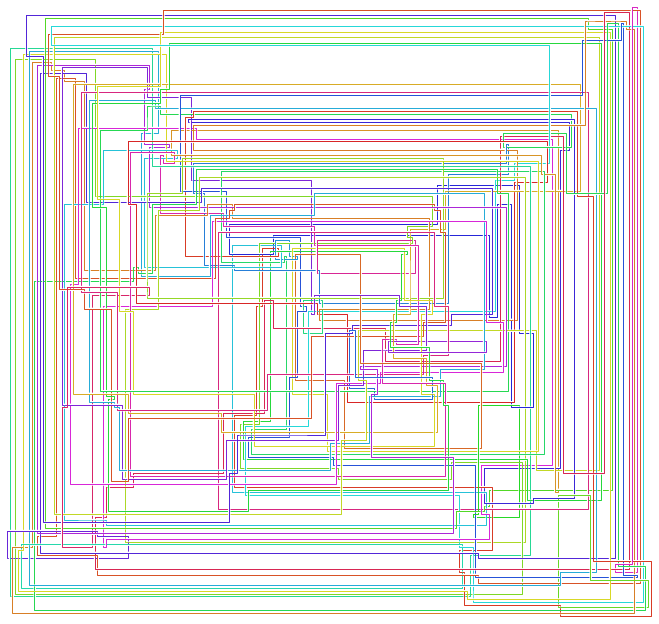}
\caption[caption]{$(1,\myLangle 3 + 2\sqrt{-1}\myRangle )$. \cite{DOR:linkDiag}.}
\end{center}
\end{figure}

\clearpage

\section{$d=2$}

Missing: $(2,\myLangle 3+\sqrt{-2}\myRangle)$.

\begin{figure}[h]
\begin{center}
\includegraphics[scale=0.5]{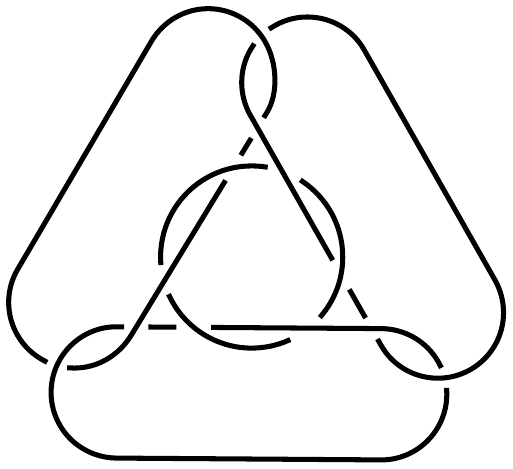}
\caption[caption]{$(2,\myLangle 1+\sqrt{-2}\myRangle)$. \cite[Example~6.8.10]{thurston:notes} \cite[Fig.~5a]{bakerReid:prinCong}.}
\end{center}
\end{figure}

\begin{figure}[h]
\begin{center}
\includegraphics[scale=0.65]{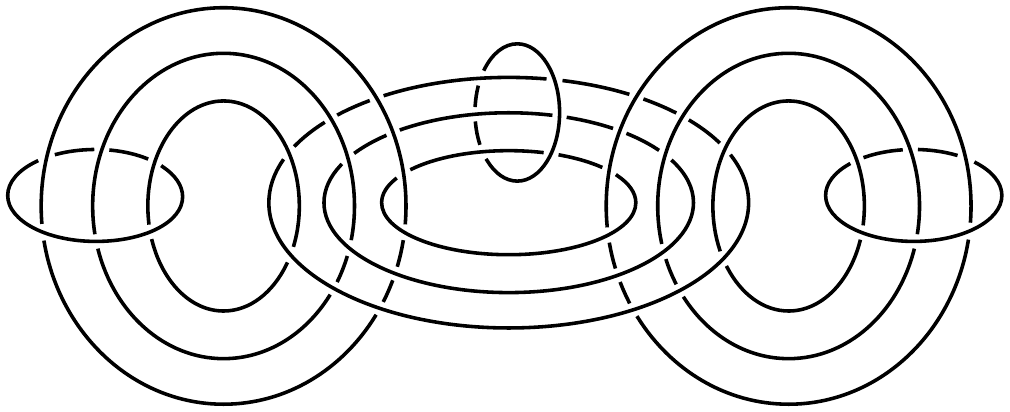}
\caption[caption]{$(2,\myLangle 2\myRangle)$. \cite[Fig.~I-13]{baker:thesis} \cite[Fig.~1b]{bakerReid:prinCong}.}
\end{center}
\end{figure}

\begin{figure}[h]
\begin{center}
\includegraphics[scale=0.45]{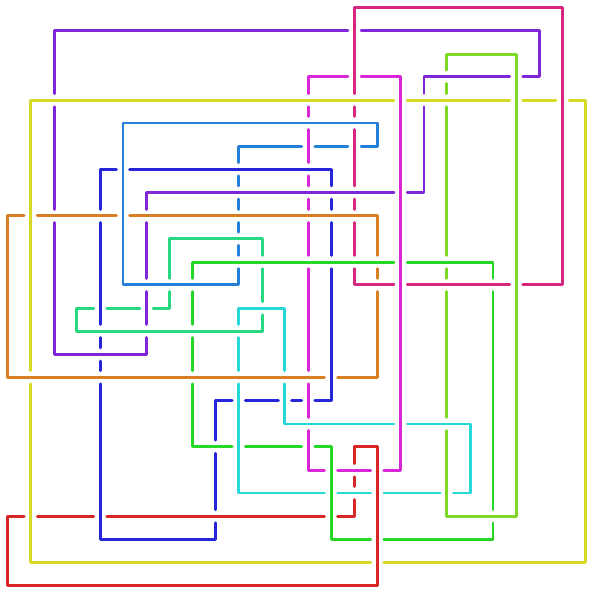}
\caption[caption]{$(2,\myLangle 2+ \sqrt{-2}\myRangle )$. \cite{DOR:linkDiag}.}
\end{center}
\end{figure}

\begin{figure}[h]
\begin{center}
\includegraphics[scale=0.2]{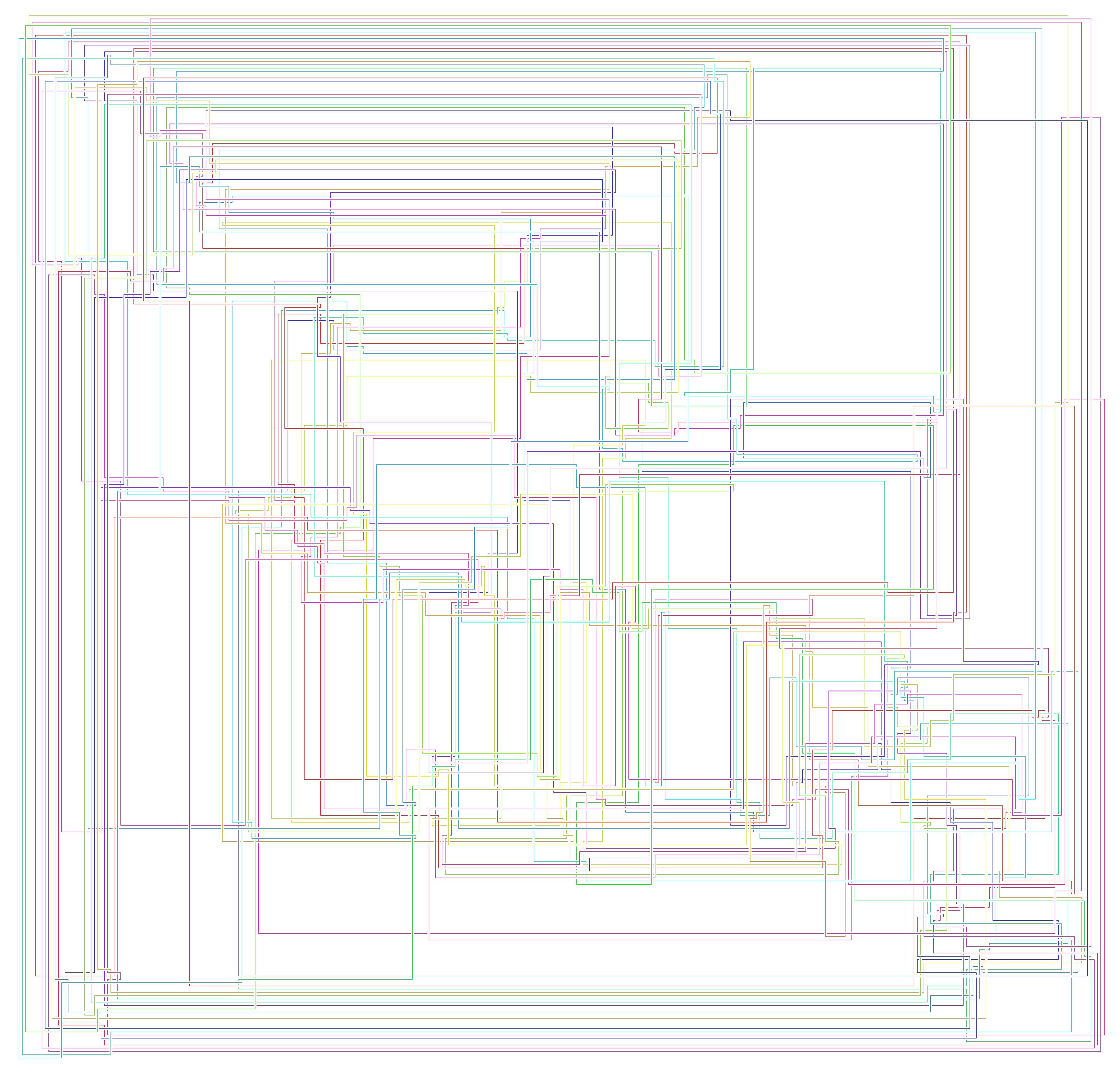}
\caption[caption]{$(2,\myLangle 1 + 2\sqrt{-2}\myRangle )$. \cite{DOR:linkDiag}.}
\end{center}
\end{figure}

\clearpage

\section{$d=3$}

Missing: $(3,\myLangle 4+\sqrt{-3}\myRangle)$, $(3,\myLangle \frac{9+\sqrt{-2}}{2}\myRangle)$.

\begin{figure}[h]
\centering
\begin{minipage}[t]{0.48\textwidth}
\begin{center}
\includegraphics[scale=1.1]{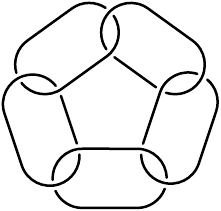}
\caption[caption]{$(3,\myLangle 2\myRangle )$.\\ \cite[Fig.~I-9]{baker:thesis}\\ \cite[Fig.~2]{dunfield:virtHaken}.}
\end{center}
\end{minipage} %
\begin{minipage}[t]{0.48\textwidth}
\begin{center}
\includegraphics[scale=1]{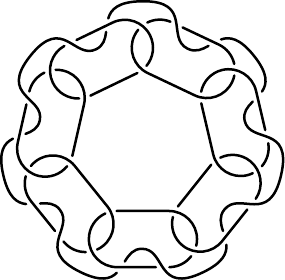}
\caption[caption]{$(3,\myLangle \frac{5+\sqrt{-3}}{2}\myRangle )$.\\ \cite{thurston:howToSee}. Also see \cite{ago:thuCongLink}.}
\end{center}
\end{minipage}
\end{figure}

\begin{figure}[h]
\begin{center}
\includegraphics[scale=0.85]{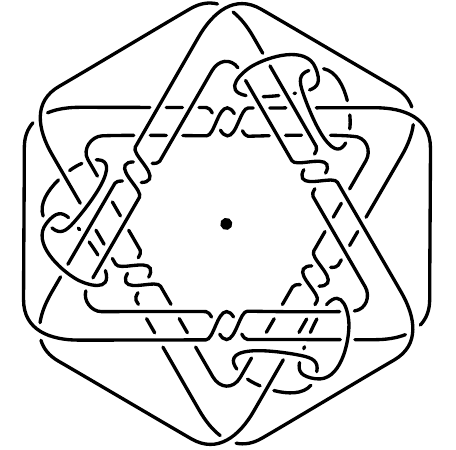}
\caption{$(3,\myLangle 3\myRangle )$. \cite[Fig.~1.7]{goerner:thesis} \label{fig:3plus0sqrtminus3}.}
\end{center}
\end{figure}

\begin{figure}[h]
\begin{center}
\includegraphics[scale=0.33]{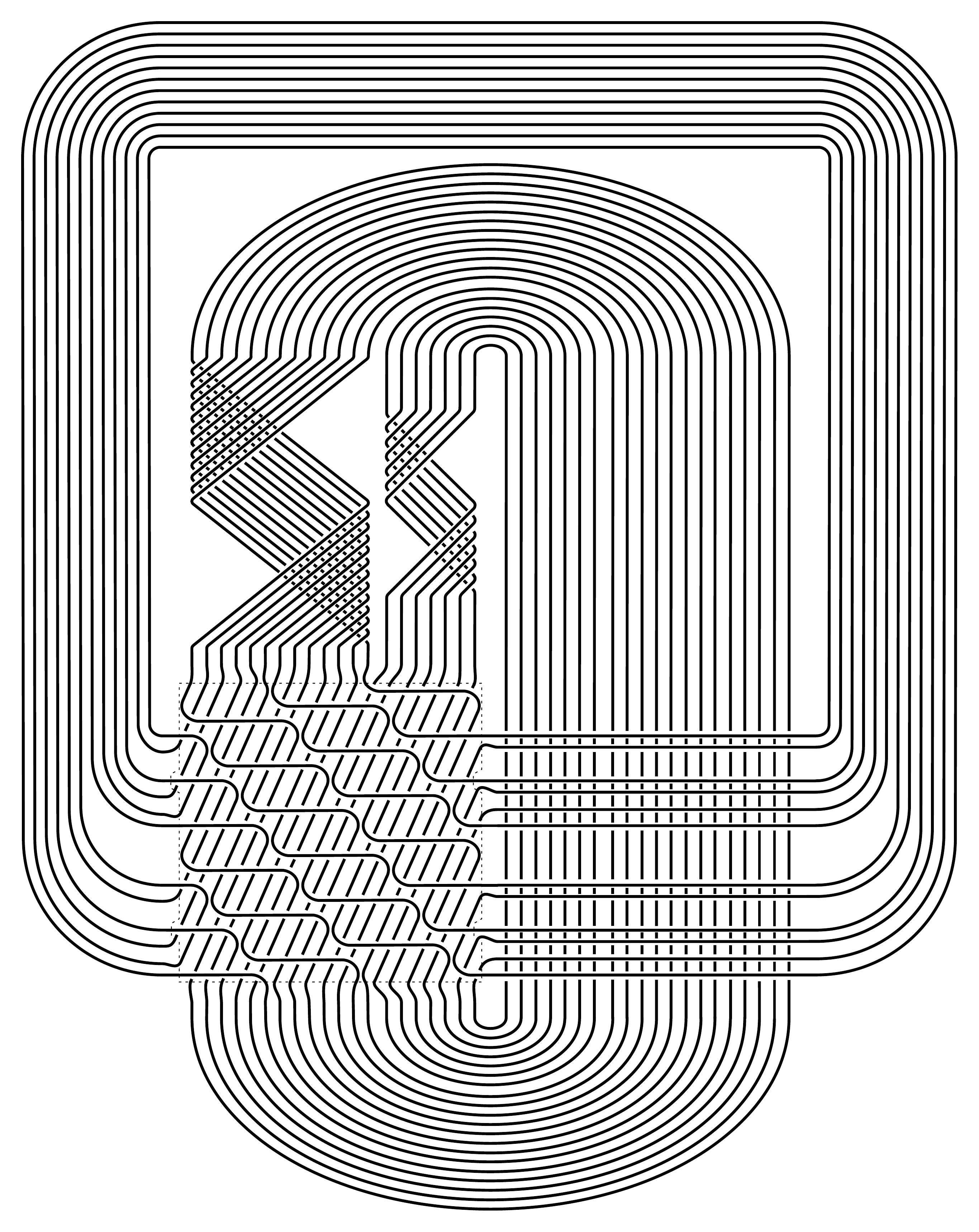}
\end{center}
\caption{$(3,\myLangle 3+\sqrt{-3}\myRangle )$. \cite[Fig.~1.27]{goerner:thesis} \label{fig:3plus1sqrtminus3}.}
\end{figure}

\begin{figure}[h]
\begin{center}
\includegraphics[scale=1.4]{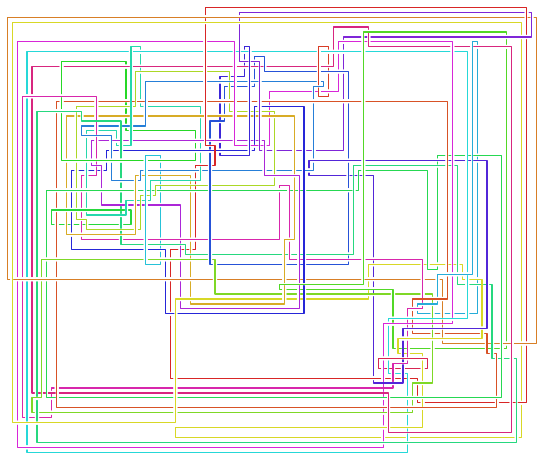}
\caption[caption]{$(3,\myLangle \frac{7 + \sqrt{-3}}{2}\myRangle )$. \cite{DOR:linkDiag}.}
\end{center}
\end{figure}

\clearpage

\section{$d=5$}

\begin{figure}[h]
\begin{center}
\includegraphics[scale=0.9]{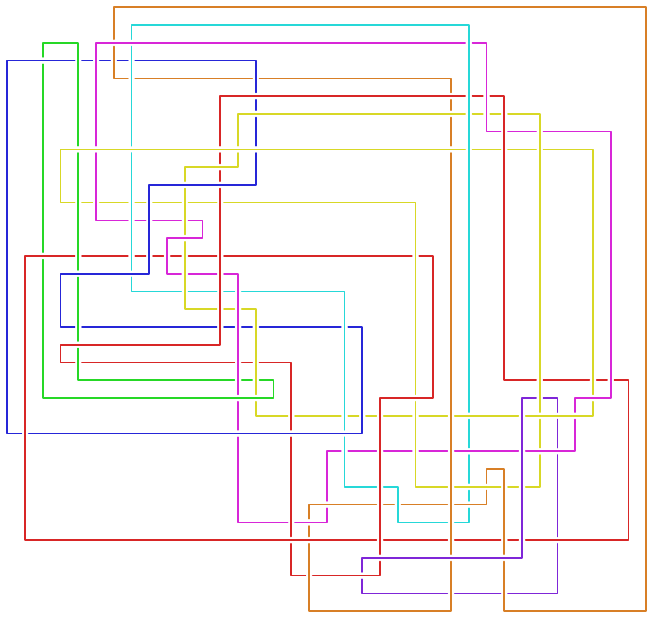}
\caption[caption]{$(5,\myLangle 3, 1 + \sqrt{-5}\myRangle )$. \cite{DOR:linkDiag}.}
\end{center}
\end{figure}

\clearpage

\section{$d=7$}

Missing: $(7,\myLangle 2+\sqrt{-7}\myRangle)$, $(7,\myLangle \frac{7+\sqrt{-7}}{2}\myRangle)$, $(7,\myLangle \frac{1+3\sqrt{-7}}{2}\myRangle)$.

\begin{figure}[h]
\centering
\begin{minipage}[t]{0.49\textwidth}
\begin{center}
\includegraphics[scale=0.5]{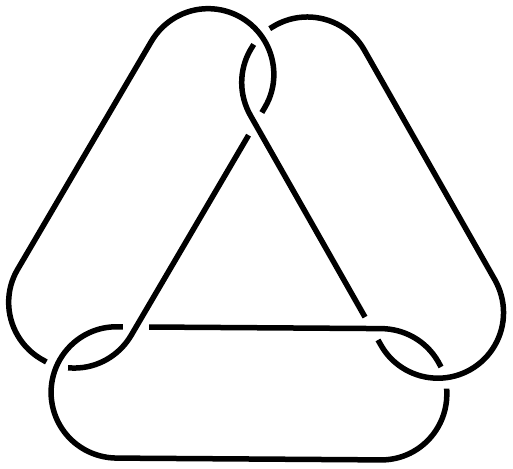}
\caption[caption]{$(7,\myLangle \frac{1+\sqrt{-7}}{2}\myRangle )$.\\ \cite[Fig~1(a)]{grunewaldSchwermer}\\ \cite[Fig.~3]{bakerReid:prinCong}.}
\end{center}
\end{minipage}
\begin{minipage}[t]{0.49\textwidth}
\begin{center}
\includegraphics[scale=0.5]{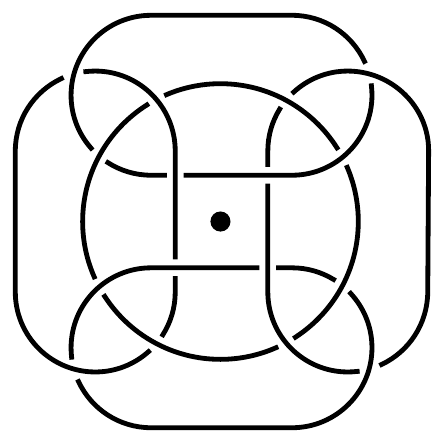}
\caption[caption]{$(7,\myLangle \frac{3+\sqrt{-7}}{2}\myRangle )$.\\ By second author.}
\end{center}
\end{minipage}
\end{figure}

\begin{figure}[h]
\begin{center}
\includegraphics[scale=0.5]{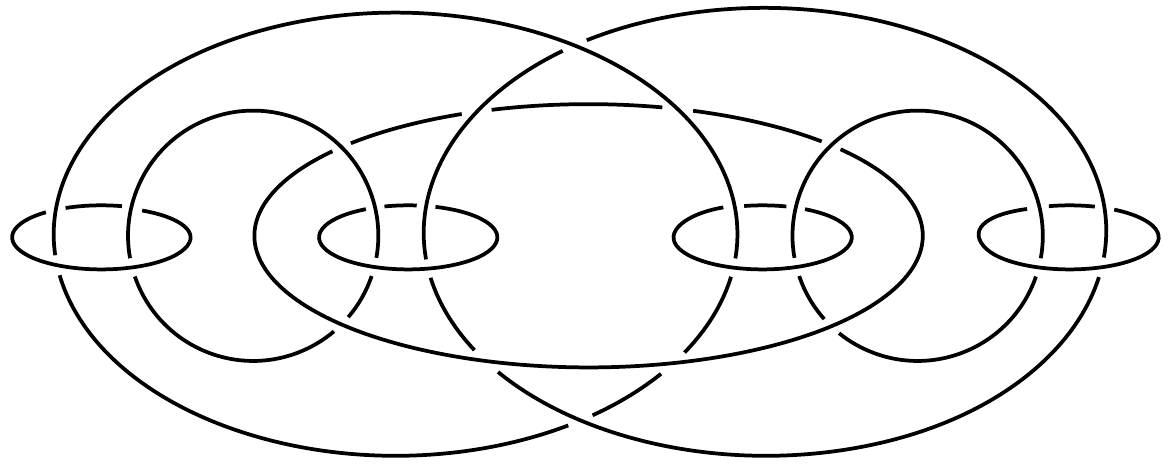}
\caption{$(7,\myLangle 2\myRangle )$. \cite[Fig.~I-5]{baker:thesis} \cite[Fig.~1d]{bakerReid:prinCong}.}
\end{center}
\end{figure}

\begin{figure}[h]
\begin{center}
\includegraphics[scale=0.8]{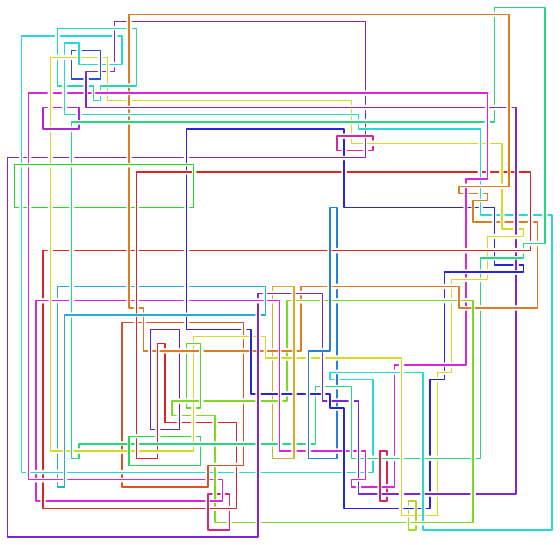}
\caption[caption]{$(7,\myLangle \sqrt{-7}\myRangle )$. \cite{DOR:linkDiag}.}
\end{center}
\end{figure}

\begin{figure}[h]
\begin{center}
\includegraphics[scale=1.2]{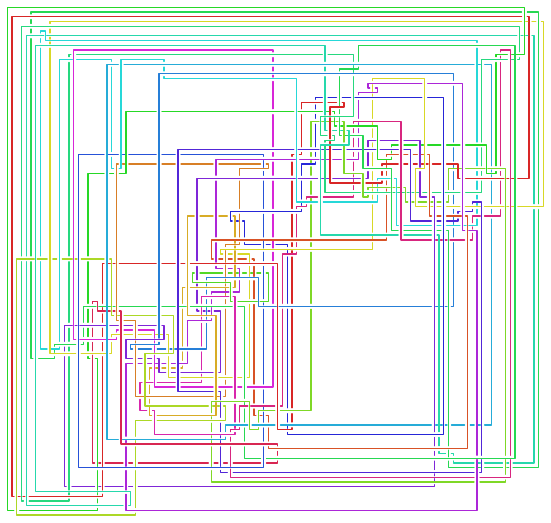}
\caption[caption]{$(7,\myLangle \frac{5+\sqrt{-7}}{2}\myRangle )$. \cite{DOR:linkDiag}.}
\end{center}
\end{figure}

\begin{figure}[h]
\begin{center}
\includegraphics[scale=0.7]{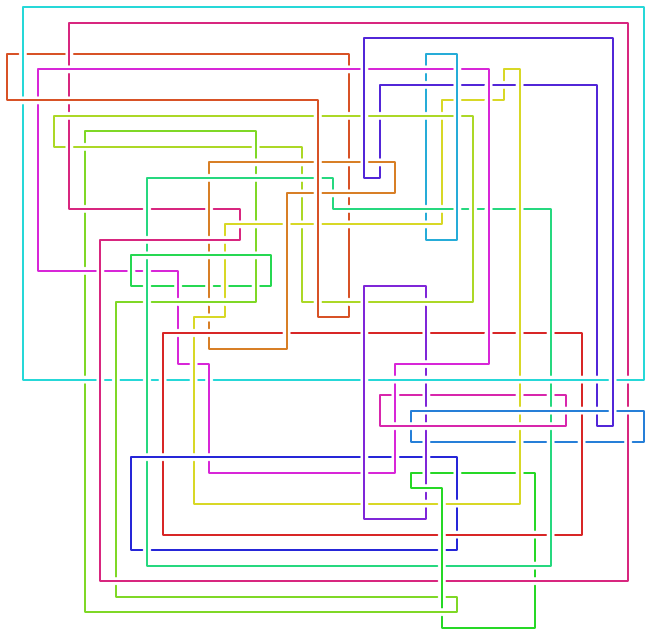}
\caption[caption]{$(7,\myLangle 1+\sqrt{-7}\myRangle )$. \cite{DOR:linkDiag}.}
\end{center}
\end{figure}

\clearpage

\section{$d=11$}

Missing: $(11,\myLangle \frac{5+\sqrt{-11}}{2}\myRangle)$.

\begin{figure}[h]
\begin{center}
\includegraphics[scale=0.5]{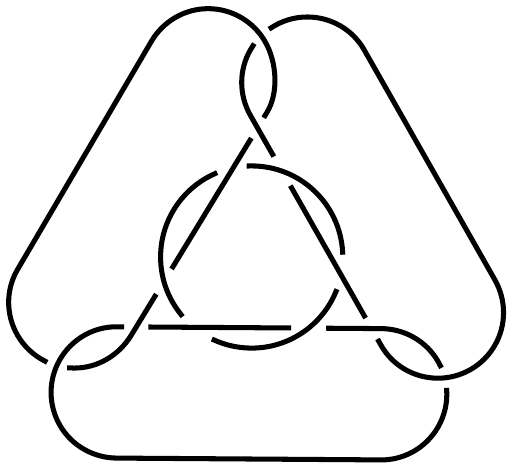}
\caption{$(11,\myLangle \frac{1+\sqrt{-11}}{2}\myRangle )$. \cite[Fig.~9]{hatcher:someLinks} \cite[Fig.~5b]{bakerReid:prinCong}.}
\end{center}
\end{figure}

\begin{figure}[h]
\begin{center}
\includegraphics[scale=1.05]{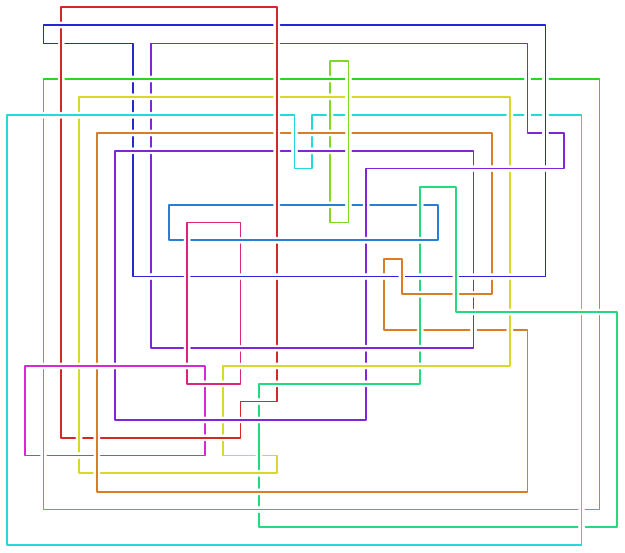}
\caption[caption]{$(11,\myLangle \frac{3+\sqrt{-11}}{2}\myRangle )$. \cite{DOR:linkDiag}.}
\end{center}
\end{figure}

\clearpage

\section{$d=15$}

\begin{figure}[h]
\begin{center}
\includegraphics[scale=0.8]{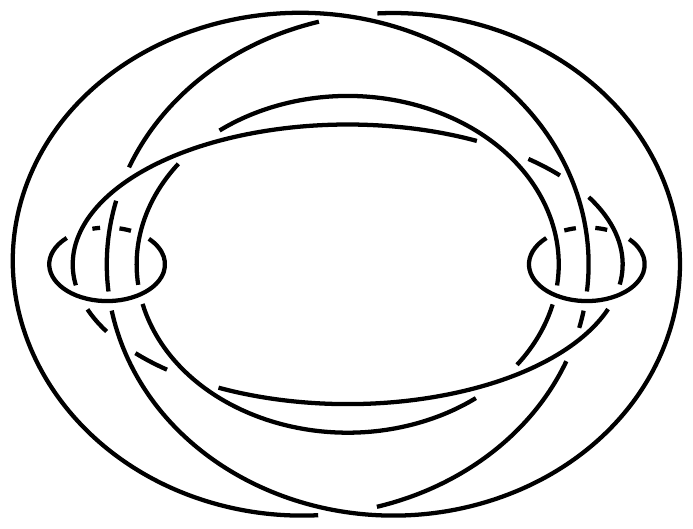}
\caption{$(15,\myLangle 2, \frac{1+\sqrt{-15}}{2}\myRangle )$. \cite[Fig.~1]{baker:top90} \cite[Fig.~4a]{bakerReid:prinCong}.}
\end{center}
\end{figure}

\begin{figure}[h]
\begin{center}
\includegraphics[scale=0.55]{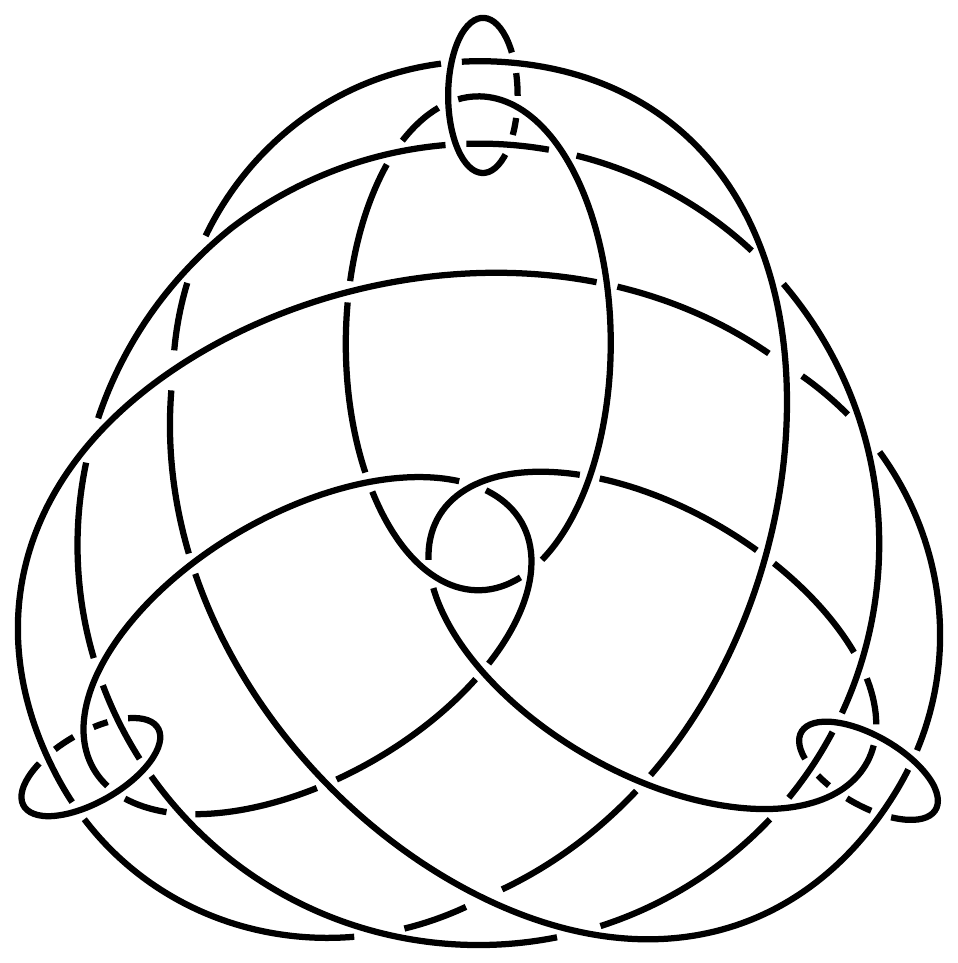}
\caption{$(15,\myLangle 3, \frac{3+\sqrt{-15}}{2}\myRangle )$. Symmetrized by second author from diagram created by \cite{DOR:linkDiag}.}
\end{center}
\end{figure}

\begin{figure}[h]
\begin{center}
\includegraphics[scale=0.6]{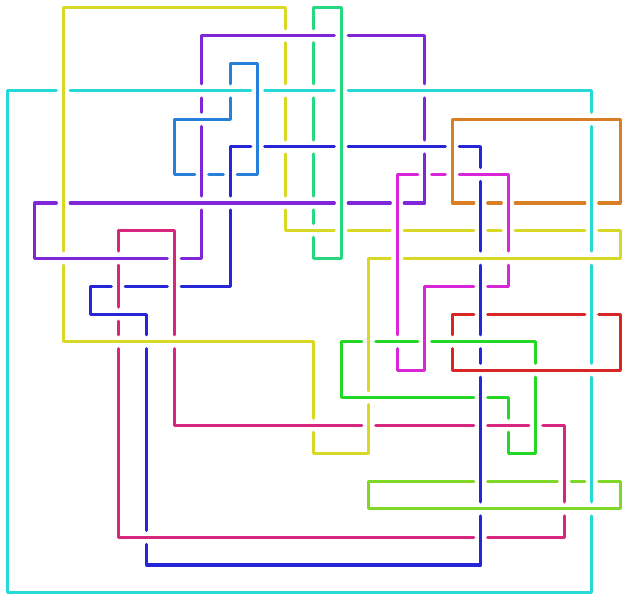}
\caption[caption]{$(15,\myLangle \frac{1+\sqrt{-15}}{2}\myRangle )$. \cite{DOR:linkDiag}.}
\end{center}
\end{figure}

\begin{figure}[h]
\begin{center}
\includegraphics[scale=1.05]{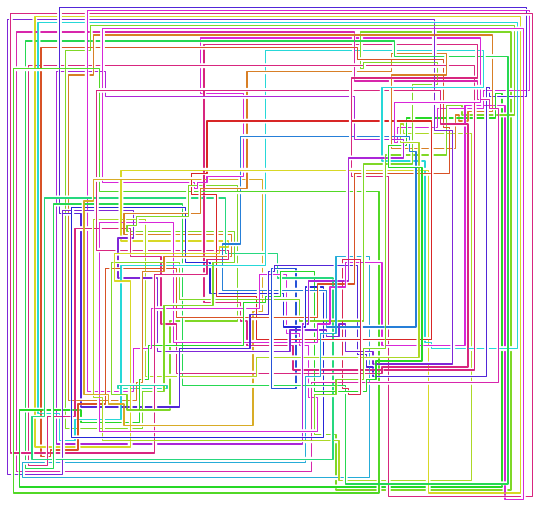}
\caption[caption]{$(15,\myLangle 5, \frac{5+\sqrt{-15}}{2}\myRangle )$. \cite{DOR:linkDiag}.}
\end{center}
\end{figure}

\begin{figure}[h]
\begin{center}
\includegraphics[scale=0.6]{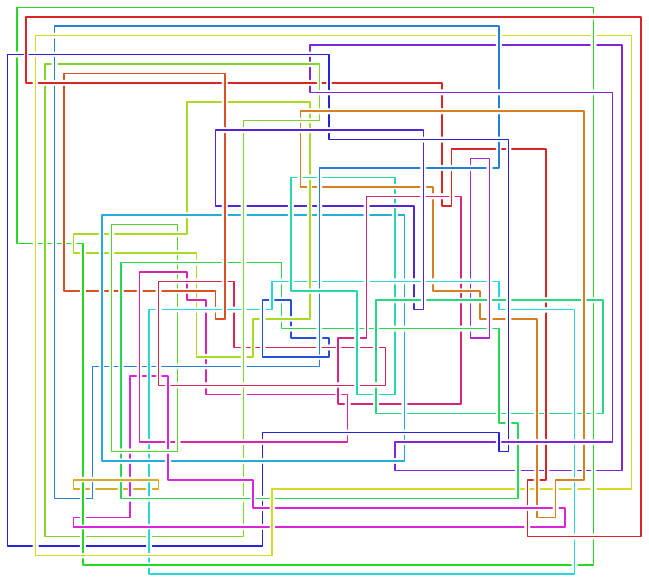}
\caption[caption]{$(15,\myLangle \frac{3+\sqrt{-15}}{2}\myRangle )$. \cite{DOR:linkDiag}.}
\end{center}
\end{figure}

\clearpage

\section{$d=19$}

\begin{figure}[h]
\begin{center}
\includegraphics[scale=0.9]{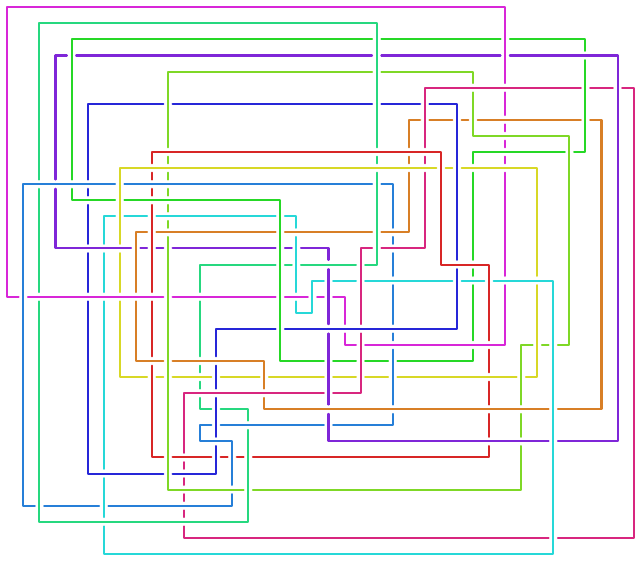}
\caption[caption]{$(19,\myLangle \frac{1+\sqrt{-19}}{2}\myRangle )$. \cite{DOR:linkDiag}.}
\end{center}
\end{figure}

\clearpage

\section{$d=23$}

\begin{figure}[h]
\begin{center}
\includegraphics[scale=0.5]{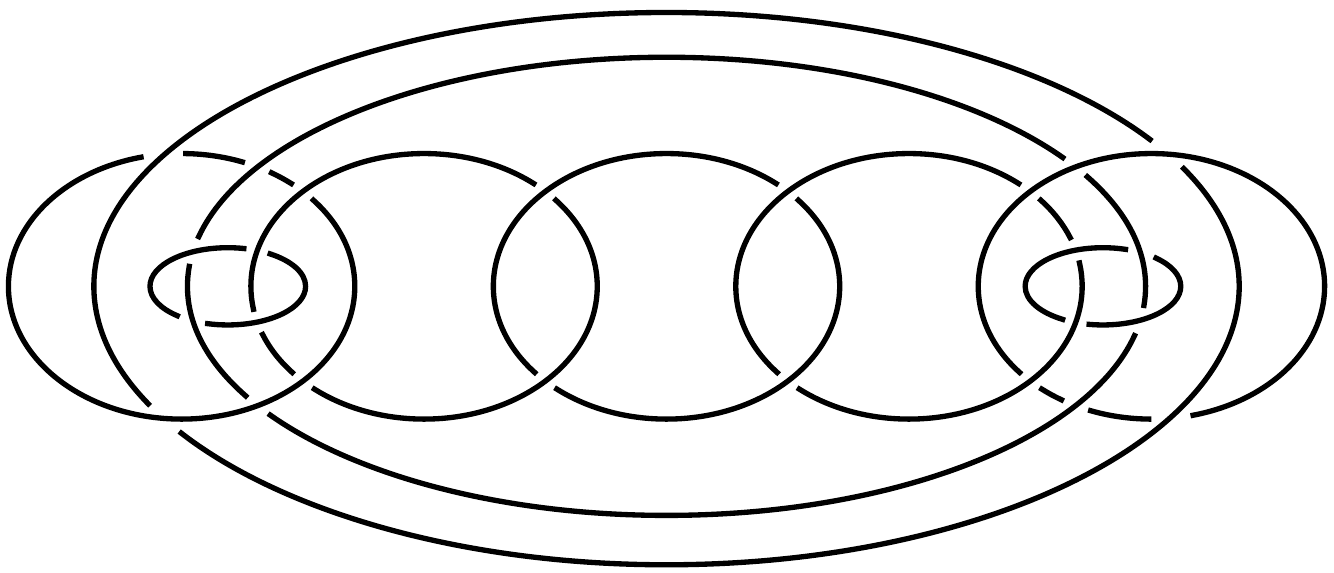}
\caption{$(23,\myLangle 2, \frac{1+\sqrt{-23}}{2}\myRangle )$. \cite[Fig.~2]{baker:top90} \cite[Fig.~4b]{bakerReid:prinCong}.}
\end{center}
\end{figure}

\begin{figure}[h]
\begin{center}
\includegraphics[scale=1.0]{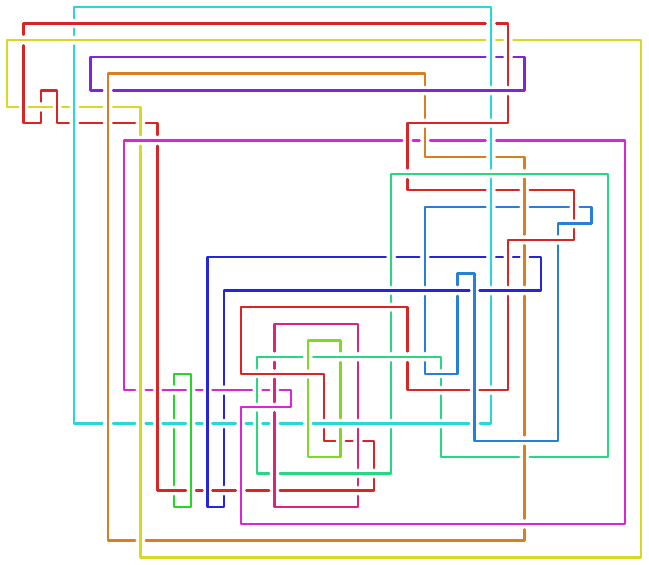}
\caption[caption]{$(23,\myLangle 3, \frac{1+\sqrt{-23}}{2}\myRangle )$. \cite{DOR:linkDiag}.}
\end{center}
\end{figure}

\begin{figure}[h]
\begin{center}
\includegraphics[scale=1.4]{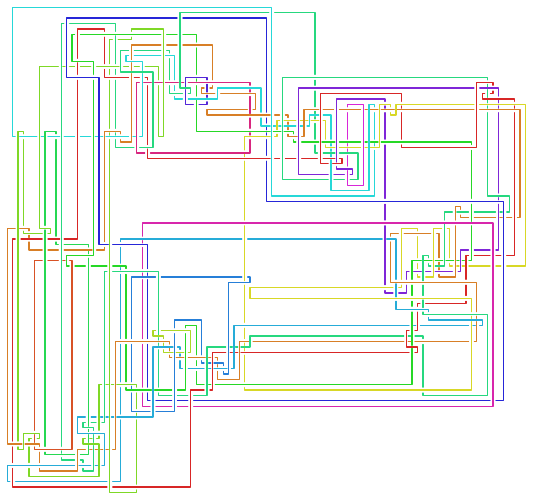}
\caption[caption]{$(23,\myLangle 4, \frac{3+\sqrt{-23}}{2}\myRangle )$. \cite{DOR:linkDiag}.}
\end{center}
\end{figure}

\clearpage

\section{$d=31$}

Missing: $(31,\myLangle 5,\frac{3+\sqrt{-31}}{2}\myRangle)$.

\begin{figure}[h]
\begin{center}
\includegraphics[scale=1]{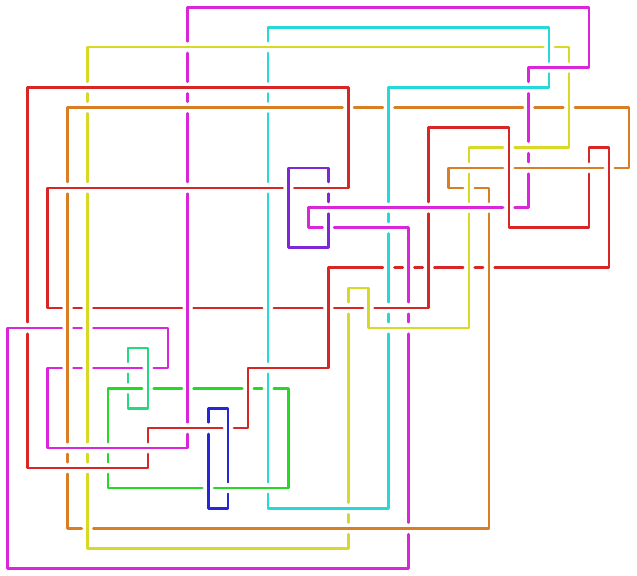}
\caption[caption]{$(31,\myLangle 2, \frac{1+\sqrt{-31}}{2}\myRangle )$. \cite{DOR:linkDiag}.}
\end{center}
\end{figure}

\begin{figure}[h]
\begin{center}
\includegraphics[scale=1.7]{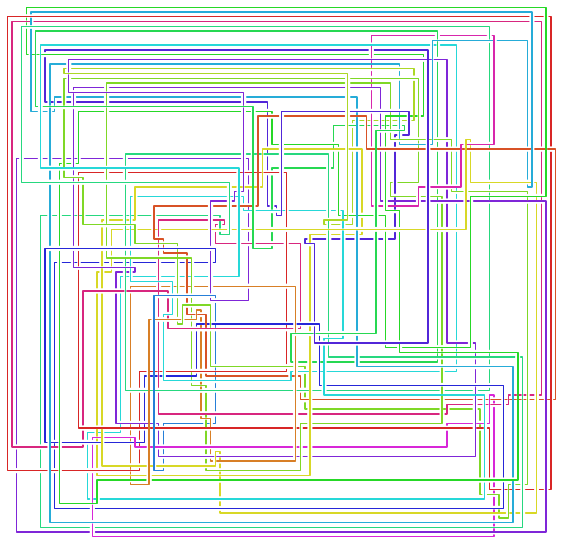}
\caption[caption]{$(31,\myLangle 4, \frac{1+\sqrt{-31}}{2}\myRangle )$. \cite{DOR:linkDiag}.}
\end{center}
\end{figure}

\clearpage

\section{$d=47$}

Missing: $(47,\myLangle 4, \frac{1+\sqrt{-47}}{2}\myRangle)$.

\begin{figure}[h]
\begin{center}
\includegraphics[scale=0.65]{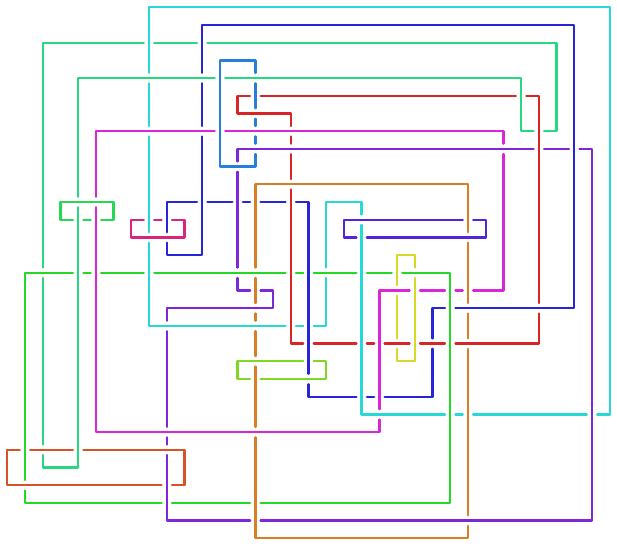}
\caption[caption]{$(47,\myLangle 2, \frac{1+\sqrt{-47}}{2}\myRangle )$. \cite{DOR:linkDiag}.}
\end{center}
\end{figure}

\begin{figure}[h]
\begin{center}
\includegraphics[scale=1.05]{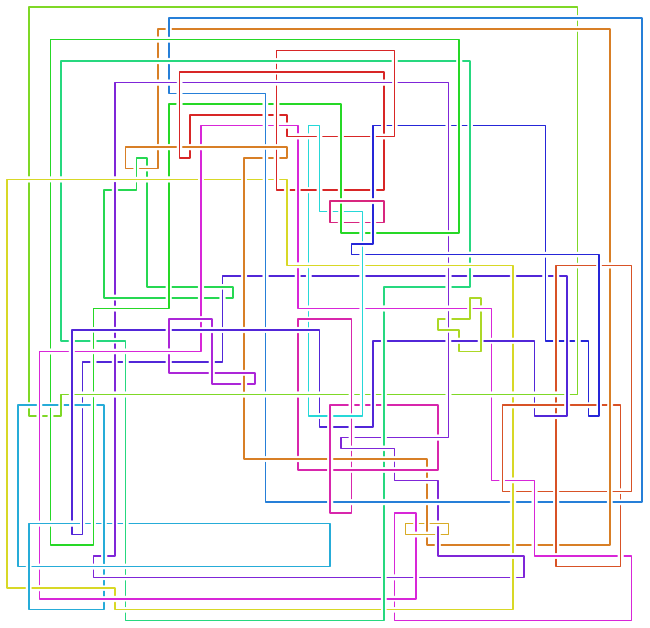}
\caption[caption]{$(47,\myLangle 3, \frac{1+\sqrt{-47}}{2}\myRangle )$. \cite{DOR:linkDiag}.}
\end{center}
\end{figure}

\clearpage

\section{$d=71$}

\begin{figure}[h]
\begin{center}
\includegraphics[scale=1.1]{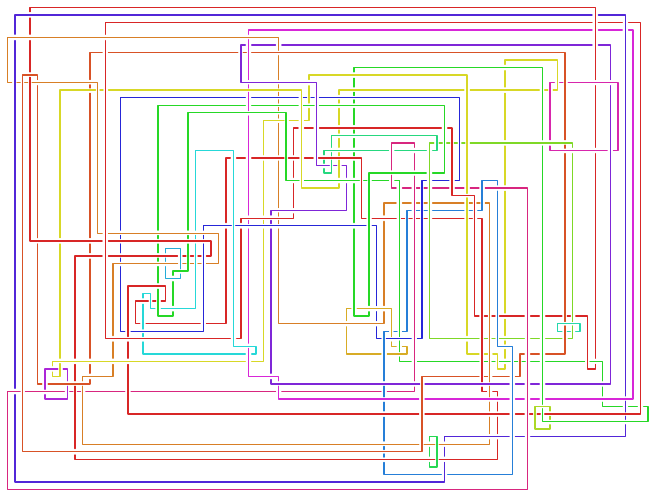}
\caption[caption]{$(71,\myLangle 2, \frac{1+\sqrt{-71}}{2}\myRangle )$. \cite{DOR:linkDiag}.}
\end{center}
\end{figure}

\clearpage

\section{Cubical}

\begin{figure}[h]
\begin{center}
\includegraphics{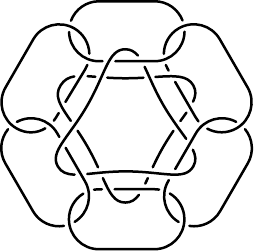}
\caption{$\MunivPrinCong{\{4,3,6\}}{\frac{3+\sqrt{-3}}{2}}$. By second author.\label{fig:cubical}}
\end{center}
\end{figure}

The link in Figure~\ref{fig:cubical} is a regular tessellation link (as defined in \cite{goerner:regTessLinkComps}) but not a principal congruence link. Similar to a principal congruence manifold, a regular tessellation manifold is a regular covering space of an orbifold. In case of a regular tessellation manifold, this orbifold is the quotient of $\H^3$ by a Coxeter group instead of a Bianchi group $\PSL(2,O_d)$. In the case of the link shown in Figure~\ref{fig:cubical}, the complement admits a tessellation into 6 regular ideal cubes such that the symmetries of the complement act transitively on flags. The complement is also arithmetic as it is also a non-regular cover of the Bianchi orbifold $Q_3=\H^3/\PSL(2,O_3)$. We refer the reader to \cite{goerner:regTessLinkComps} for further details.

\bibliographystyle{hepMatthias}
\bibliography{linkFigures}

\end{document}